\documentclass[11pt]{article}

\usepackage{amsfonts}
\usepackage{a4}
\usepackage{pictex}
\usepackage{amssymb}

\def\<{\langle}
\def\>{\rangle}
\def\Z{\mathbb Z}

\def\ar{\rightarrow}

\def\qed{{\hfill $\Box$}}

\def\ar{\rightarrow}

\def\NN{{\mathcal N}}
\def\VV{{\mathcal V}}
\def\CC{{\mathcal C}}

\def\proof{\noindent{\bf Proof.\ }}

\makeatletter
\newcommand{\fpd}{\mathop{\operator@font \hbox{\Large$\ast$}}}
\makeatother

\newtheorem{lemma}{Lemma}
\newtheorem{proposition}{Proposition}
\newtheorem{corollary}{Corollary}
\newtheorem{theorem}{Theorem}

\title{On the profinite topology of right-angled Artin groups}
\author{V. Metaftsis\\Department of Mathematics\\University of the Aegean\\Karlovassi\\832 00 Samos\\ 
Greece\\{\it email: vmet@aegean.gr} \and E. Raptis\\Department of Mathematics\\University of Athens\\ Panepistimiopolis\\157 84 Athens\\Greece\\{\it email: eraptis@math.uoa.gr}}  
\date{}

\begin{document}

\maketitle
\begin{abstract} 
In the present work, we give necessary and sufficient conditions on the graph
of a right-angled Artin group that determine whether the group is subgroup
separable or not. Moreover, we investigate the profinite topology of $F_2\times
F_2$ and we show that the profinite topology of the above group is strongly connected with the profinite
topology of $F_2$. 
\end{abstract}

\section{Introduction}

Subgroup separability is an extremely powerful property of groups with many 
topological implications. As shown by Thurston, subgroup separability allows certain 
immersions to lift to an embedding in a finite cover. Scott in \cite{scott} showed that 
subgroup separability is inherited by
subgroups and finite extensions. Although free products of
subgroup separable groups are subgroup separable, the same is not true for direct
products. This is one of the motivations for the present work.

On the other hand, although right-angled Artin groups are known for some time, (see \cite{green,droms})
they recently attracted special attention. Bestvina
and Brady \cite{bb} used the kernels of their epimorphisms to $\Z$ to construct examples 
of groups with strange finiteness properties amongst other things. 

Charney and Davis \cite{cd} and Meier and VanWyk \cite{mv} constructed, from the
graph $G$, a cubical complex (CW-complex) and they proved that it is in fact
the Eilenberg-MacLane space of $G$. Hsu and Wise \cite{hw2} showed that
$G$ is a coxeter group and Papadima and Suciu \cite{ps} calculated various
algebraic invariants for $G$ including the lower central series quotients. Also,
Meier, Meinert and VanWyk \cite{mmv} determined their geometric invariants introduced 
by Bieri, Neummann and Strebel.

In the present paper we study the profinite topology of a right-angled Artin group $G$
and we show that one can decide if $G$ is subgroup separable or not by just examining its 
graph. Moreover, we show that the only obstructions for $G$ to be subgroup separable
are the two well known examples of non-subgroup separable groups $F_2\times F_2$ and
$L$ (see \cite{ls} and \cite{nw} respectively). This was the motivation to
study the profinite topology of $F_2\times F_2$ and of the $BKS$ group (see \cite{bks})
which is responsible for the non-subgroup separability of $L$. It turned out that, 
for the $F_2\times F_2$ case, the problem of determining all closed
subgroups in its profinite topology is equivalent to determining the residual finiteness
of every finitely presented group. Nonetheless, the positive result is that
all finitely presented subgroups of $F_2\times F_2$ are closed in the profinite
topology of $F_2\times F_2$.  In fact
our results show that the profinite topologies of $F_2\times F_2$, is
strongly connected with that of $F_2$.

\section{Notation and definitions}

In this section we establish notation and we review some basic definitions and results.

By a graph $X$ we mean a finite simplicial graph with vertex set $VX$ and edge set
$EX$. The full subgraph $Y$ of $X$ is a graph whose vertex set
is a subset of $VX$, two vertices in $Y$ being adjacent in $Y$ 
if and only if they are adjacent in $X$. So the full subgraphs of
a graph $X$ are uniquely determined by their vertex sets. In the sequel,
by a subgraph $Y$ of a graph $X$, we mean the full subgraph of $X$,
defined by $VY$. 

If $X$ is a connected graph, we make $VX$ a metric space by assuming that
the length of each edge is 1. So, a full subgraph $Y$ of $X$ is a path
of length $n$, if $Y$ is the graph

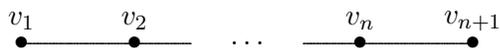
\begin{center}
\hskip 2cm
\beginpicture
\setcoordinatesystem units <1.5cm, 1cm>
\put{$\bullet$} at 3.5 0
\put{$v_1$} at 3.5 .3
\put{$\bullet$} at 4.5 0
\put{$v_2$} at 4.5 .3
\put{$v_n$} at 6.5 .3
\put{$v_{n+1}$} at 7.5 .3 
\put{$\bullet$} at 6.5 0
\put{$\bullet$} at 7.5 0
\plot 3.5 0 4.5 0 /
\plot 6.5 0 7.5 0 /
\plot 4.5 0 5 0 /
\plot 6 0 6.5 0 /
\put{$\ldots$} at 5.5 0
\endpicture 
\end{center} 
If $v_1=v_{n+1}$ we say that $Y$ is a {\it closed path\/} of length $n$.
By a {\it square\/} we mean a closed path of length 4.

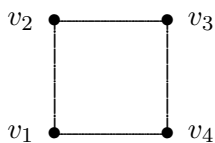
\begin{center}
\hskip 2cm
\beginpicture
\setcoordinatesystem units <1.5cm, 1.5cm>
\put{$\bullet$} at 3.5 0
\put{$\bullet$} at 3.5 1
\put{$\bullet$} at 4.5 0
\put{$\bullet$} at 4.5 1
\put{$v_1$} at 3.2 0
\put{$v_2$} at 3.2 1
\put{$v_4$} at 4.8 0
\put{$v_3$} at 4.8 1
\plot 3.5 0 3.5 1 /
\plot 3.5 1 4.5 1 /
\plot 4.5 1 4.5 0 /
\plot 4.5 0 3.5 0 /
\endpicture
\end{center}

A {\it homeomorphism\/} between graphs is a simplicial function that
is one-to-one on both vertices and edges and preserves adjacency.

Let $X$ be a finite simplicial graph. The {\it graph group\/} or 
the {\it right-angled Artin group\/} $G(X)$ (or $G$ for simplicity)
is given by the presentation with a generator $g_i$ for every
vertex $v_i$ of $X$ and a defining relation $[g_i,g_j]=1$ for
each edge between vertices $v_i$ and $v_j$ in $X$. 

Let $X$ be a graph and $G(X)$ its right-angled Artin group. Let
also $Y$ be a subgraph of $X$. Then we can also define the right-angled 
Artin group of $Y$, $G(Y)$ and it is obvious that there is a natural
embedding $G(Y) \ar G(X)$. Hence, without loss of generality, from now
on we will consider $G(Y)$ as a subgroup of $G(X)$.

The {\it profinite topology\/} of $G$ is the topology whose base of closed sets
consists of the finite index normal subgroups of $G$. Given the profinite
topology, $G$ is of course a topological group (the group operations
are continuous) and it is residually finite if and only if it is
Hausdorff (the trivial subgroup is closed with respect to the profinite
topology). A subset $H\subseteq G$ is {\it separable} in $G$ if it
is closed in the profinite topology of $G$. One can easily show that
if $K<H<G$ with $|H:K|< \infty$ then if $K$ is closed in the profinite
topology of $G$, so is $H$. 

A group $G$ is called {\it cyclic subgroup separable\/} (or $\pi_c$) if every cyclic
subgroup of $G$ is closed in the profinite topology of $G$. A group $G$ is called
{\it subgroup separable\/} (or {\it LERF\/}) if all its finitely generated subgroups are
separable. Moreover, every subgroup of a subgroup separable group
is subgroup separable \cite{scott}. Subgroup separability is a ``rare'' property of
groups. A list of known subgroup separable groups can be found in \cite{???}.

On the other hand, non-subgroup separability is also difficult to prove. We
give here two well known examples of non-subgroup separable groups that play
a major r\^ole in the sequel.

By $L$ we denote the group with presentation
$$L=\< a,b,c,d \mid [a,b]=[b,c]=[c,d]=1\>.$$
$L$ was shown to be non-subgroup separable by Niblo and Wise in \cite{nw}. In fact,
it was shown that $L$ contains a subgroup isomorphic to an index two subgroup
of the famous example of Burns, Karrass and Solitar \cite{bks}, the group with
presentation
$$BKS=\< t, a,b \mid [a,b]=1, tat^{-1}=b\>.$$

The second example is older. If $F_2$ denotes the free group of rank two
then the group $F_2\times F_2$ was shown by Michailova (see \cite{ls})
to have non-solvable generalized word problem. Consequently, 
$F_2\times F_2$ is not subgroup separable.

Finally, let $f:G\ar G$ be an automorphism of $G$. Then Fix$(f)=\{ g\in G \mid f(g)=g\}$.
Obviously, Fix$(G)$ is a subgroup of $G$.

\section{Subgroup separability}

All right-angled Artin groups are residually finite by the work of Green
\cite{green} and linear by the work of Humphries \cite{humphries}. In fact
they are $\Z$-linear by the work of Hsu and Wise \cite{hw2} and Brown \cite{brown}. 

\begin{theorem}
All polycyclic subgroups of a right-angled Artin group $G$ are closed in
the profinite topology of $G$. In particular, $G$ is cyclic subgroup separable.
\end{theorem}

\proof  Let $G$ be a right-angled Artin group. Then $G$ is linear and in fact,
$G$ is a subgroup of GL$(n,\Z)$. Hence, by \cite[Corollary 1, page 26]{segal},
every soluble subgroup of GL$(n,\Z)$ is polycyclic and so is every soluble
subgroup of $G$. But all polycyclic subgroups of GL$(n,\Z)$ are closed
in the profinite topology of GL$(n,\Z)$ (see \cite[Theorem 5, page 61]{segal}). 
Therefore, every polycyclic subgroup
of $G$ is closed in the subspace topology of $G$ which is coarser than the
profinite topology of $G$. Consequently, every
cyclic subgroup of $G$ is closed in the profinite topology of $G$, so 
$G$ is cyclic subgroup separable. \qed

\begin{lemma}\label{path}
Let $G$ be a right-angled Artin group with graph $X$. 
If $X$ has a path of length three as a subgraph then $G$
is not subgroup separable.
\end{lemma}

\proof If $X$ has a subgraph $Y$ homeomorphic to a path of length three
then $G(Y)$ is isomorphic to $L$ and so $G$ has a subgroup isomorphic to $L$ 
and hence cannot be subgroup
separable. \qed

\begin{lemma}\label{square}
Let $G$ be a right-angled Artin group with graph $X$. If $X$ has a
subgraph $T$ which is a closed path of length four or more, 
then $G$ is not subgroup separable.
\end{lemma}

\proof If $T$ has length five or more then $T$ contains a subtree with
a path of length at least three and so $T$ and hence $X$ contain
a subgraph homeomorphic to a path of length three thus $G$ cannot
be subgroup separable by Lemma \ref{path}.  

Else, the subgraph $T$ is homeomorphic to a square with vertices
$v_a,v_b,v_c,v_d$. Then the right-angled group $G(T)$ is the subgroup
of $G$ generated by $\< a,b,c,d\>$,  with presentation
$$G(T)=\< a, b, c,d\mid [a,b]=[b,c]=[c,d]=[d,a]=1\>=\< a,c\> \times \<b,d\>,$$ 
hence $G(T)$ is isomorphic to $F_2\times F_2$ where $F_2$ is the free group of 
rank two. This last group is well known to be non subgroup separable by the
work of Michailova \cite{ls}.\qed

\begin{lemma}\label{piliko}
Let $G$ be a right angled Artin group with connected graph $X$. 
If $v_a$ is a vertex of $X$ connected to every other vertex of $X$
then $G=R\times \< a\>$, where $R$ is the right angled Artin group with
graph the full subgraph of $X$ with vertex set $VX\setminus\{ v_a\}$. 
\end{lemma}

\proof Since $v_a$ is connected to every other vertex of $X$,
we have $G=R\times \Z$ where $R$ is the subgroup
of $G$ generated by all the generator of $G$ but $a$. Obviously,
$R$ contains the relations of $G$ that do not involve $a$. So,
in graph theoretic language, $R$ involves all vertices of
$X$ but $v_a$ as well as all edges of $X$ but those that connect
vertices to $a$. Hence, $R$ is the subgroup of $G$ that
corresponds to the subgraph with vertex set $VX\setminus \{v_a\}$.
\qed

\begin{theorem}
Let $G$ be a right-angled Artin group with graph $X$. 
Then $G$ is subgroup separable if and only if $X$ does not contain 
a subgraph homeomorphic to either a square or a path of length three. 
\end{theorem}

\proof Without loss of generality we may assume that $X$ is connected.
If $X$ is disconnected we work with the connected components of
$X$. The subgroup separability of $G$ is then a consequence of 
the fact that the free product of two subgroup separable groups
is subgroup separable. 

Assume first that $X$ does not contain a subgraph homeomorphic to
either a square or a path of length three. We use induction on the 
number of vertices of $X$. 

If $X$ contains one or two vertices then $G$ is isomorphic to either
$\Z$ or $\Z^2$ and so is subgroup separable. If $X$ contains three vertices 
then there is at least one vertex,
say $v_a$, that is connected to every other vertex of $X$. Then
by Lemma \ref{piliko}, $G=A\times \Z$ where $A$ is either a free
abelian group of rank two or a free group of rank two. In both
cases $G$ is subgroup separable, in the first since it is abelian
and in the second, by the work of Allenby and Gregorac \cite{ag}.

Assume that every right-angled Artin group having a graph with $k$ 
vertices that contains no subgraph homeomorphic to either a 
square or a path of length three is subgroup separable.   

Let $Y$ be a graph with $k+1$ vertices that satisfies the hypotheses
of the theorem. Then by Lemma in \cite{droms}, there is at least
one vertex in $Y$ that is connected to every other vertex of $Y$.
So $R=M \times \Z$ where, by Lemma \ref{piliko},  
$M$ is a right-angled Artin group
that corresponds to the subgraph with vertex set $VY\setminus\{ v\}$. 
So $M$ is subgroup separable from the inductive hypothesis and so $R$ is subgroup
separable from Lemma 3 in \cite{mr2}. 

Conversely, if $G$ is subgroup separable it cannot contain a subgroup
isomorphic neither to $L$ nor to $F_2\times F_2$. Hence, its graph $X$
cannot have a subgraph homeomorphic to neither a square nor a path of
length three. \qed

We should mention here that the above theorem easily generalizes to graph groups,
that is Artin groups with each vertex associated to a
free abelian group of finite rank.

%
\section{The profinite topology of $F_2\times F_2$.}

The following Lemma is a simple generalization of Lemma 2 in \cite{mr3}. The proof is 
practically the same as of \cite[Lemma 2]{mr3} but is included here for completeness.

\begin{lemma}[\cite{mr3}]\label{palio}
Let $G$ be a group and let $H$ be a finitely generated, subgroup separable,
normal subgroup of $G$ such that $G/H'$ is subgroup separable for every
characteristic subgroup $H'$ of $H$. Let also $M$ be a finitely generated
subgroup of $G$. Then $M$ is closed in the profinite topology of $G$ if
$M\cap H$ is closed in the profinite topology of $H$.
\end{lemma}

\proof It suffices
to show that $\bigcap_{N\in\NN}MN=M$ where $\NN$ is the set of all
normal subgroups of finite index in $G$. Let $\CC$ be the set of all characteristic subgroups
of finite index in $H$. For every $H'\in\CC$ we have that $G/H'$ is subgroup separable and
that $MH'/H'$ is finitely generated, hence 
$$\bigcap_{V\in\VV}V\frac{MH'}{H'}=\frac{MH'}{H'}$$
or equivalently
$$\bigcap_{N\in\NN}\frac{NH'}{H'}\frac{MH'}{H'}=\frac{MH'}{H'}$$
where  $\VV$ is the set of all 
normal subgroups of finite index in $G/H'$. Consequently, $\bigcap\limits_{N\in\NN}MN$ 
is a subset of $MH'$ for every $H'\in\CC$.

Now, let $U=\bigcap\limits_{H'\in\CC}MH'$. Obviously $M$ is a subgroup of
$U$. So, 
$$U\cap H=(\bigcap\limits_{H'\in\CC}MH')\bigcap H=\bigcap\limits_{H'\in\CC}(M\cap H)H'.
$$ But $$\bigcap\limits_{H'\in\CC}(M\cap H)H'= \bigcap\limits_{N\in\NN'}(M\cap H)N
=M\cap H$$ since $M\cap H$ is closed in the profinite topology of $H$. In the above, 
$\NN'$ is the
set of all finite index normal subgroups of $H$. So, $U\cap H=M\cap H$.

Let $u\in U$. Then, for every $H'\in\CC$ there is an $h'\in H'$ and an $l'\in M$ such that 
$u=l'h'$. Hence, $(l')^{-1}u=h'\in H'$ and so $(l')^{-1}u\in H$. On the other hand, $M$ is a 
subgroup of $U$ and so $l'\in U$. Therefore $(l')^{-1}u\in U$.
Hence, $(l')^{-1}u\in U\cap H=M\cap H$. Thus, there is an
$l_1\in M$ such that $l^{-1}u=l_1$ which implies that $u=ll_1\in
M$. So $U\subseteq M$. But
$M\subseteq U$ and therefore
$U=M$. Since $\bigcap\limits_{N\in\NN}MN\subseteq U=M$ we have
that $\bigcap\limits_{N\in\NN}MN=M$ as required. \qed

If $C=A\times B$ then, by abusing notation, we identify $A\times \{ 1\}$ with $A$ and
$\{ 1\}\times B$ with $B$.
So we can now prove the following.

\begin{proposition}\label{tomi}
Let $C=A\times B$ where $A,B$ are subgroup separable groups. A finitely  generated subgroup 
$M$ of $C$ is closed in the profinite
topology of $C$ if and only if $M\cap A$ (or $M\cap B$) is closed in the
profinite topology of $A$ (or $B$).
\end{proposition}

\proof If $M\cap A$ is closed in the profinite topology of $A$  then $M$ 
is closed in the profinite topology of $C$, by Lemma \ref{palio}. 

Assume now that $M$ is closed in $C$. Both $A$ and $B$ are also closed in 
the profinite topology of $C$. Indeed,
if $g\in G$ with $g\not\in A$ then under the projection homomorphism 
$f: A\times B\ar B$, $f(A)=1$ but $f(g)\neq 1$. The result follows easily 
from the fact that $B$ is subgroup separable. Hence $M\cap A$ 
is closed in the profinite topology of $C$ as an intersection of closed sets. Consequently,
$M\cap A$ is closed in the subspace topology of $A$ which is coarser than the
profinite topology of $A$. Hence, $M\cap A$ is closed in the profinite
topology of $A$. The case $M\cap B$ is equivalent. \qed

Now we can use the above proposition to show a positive and a negative result.

Let $F_2'$ be an isomorphic copy of $F_2$, the free group of rank two.
The positive result is the following.

\begin{corollary}
Let $H$ be a finitely presented subgroup of $G=F_2\times F_2'$. Then $H$
is closed in the profinite topology of $G$.
\end{corollary}

\proof By the work of Baumslag and Roseblade \cite{br}, $H$ is either free
or else has a subgroup of finite index that is the product $H_1\times H_2$
with $H_1=F_2\cap H$ and $H_2=F_2'\cap H$. In the second case, 
each $H_i$, $i=1,2$ is finitely generated and so is closed in the profinite
topology of $F_2$ (and $F_2'$) so $H_1\times H_2$ is closed in the profinite
topology of $G$, by Proposition \ref{tomi}. Consequently, $H$ is closed
in the profinite topology of $G$.

In the first case, let $H$ be a free subgroup of $F_2\times F_2'$. If either $H\cap F_2$
or $H\cap F_2'$ are trivial then $H$ is closed in the profinite topology of $F_2\times F_2'$
by Proposition \ref{tomi}. If, on the other hand, $H\cap F_2\neq 1\neq H\cap F_2'$
then   $H$ contains a subgroup isomorphic to $\Z^2$, a contradiction to the hypothesis
that $H$ is free.  \qed

Now the negative result. The following construction is based on an idea of Michailova
\cite{ls}. Let $$H=\< x_1, \ldots, x_n\mid r_1,\ldots, r_m\>$$ be any finitely
presented group and let $F_n$ be the free group on abstract generators $\< x_1,\ldots, x_n\>$.
Obviously, $F_n\times F_n$ can be considered as a finite index subgroup of $F_2\times F_2'$,
so every subgroup of $F_n\times F_n$ is closed in the profinite topology of $F_n\times F_n$
if and only if it is closed in the profinite topology of $F_2\times F_2'$.

Let $L_H$ be the subgroup of $F_n\times F_n$ generated by
$$L_H=\< (x_i,x_i),\ \ i=1,\ldots, n, \ \ (1,r_j), \ \ j=1,\ldots, m\>.$$
Then $L_H\cap F_n$ is the normal closure of $\< r_j, j=1,\ldots ,m\>$
as a subgroup of $F_n$. So, by Proposition \ref{tomi},
$L_H$ is closed in the profinite topology of $F_n\times F_n$, if and only if
$L_H\cap F_n$ is closed in the profinite topology of $F_n$
or equivalently if and only if the group
$$H=\< x_1, \ldots, x_n\mid r_1,\ldots, r_m\>$$
is residually finite. So we have the following

\begin{corollary}
The problem of determining all closed finitely generated subgroups of $G=F_2\times F_2'$ 
with respect to the profinite topology is
equivalent to the problem of determining the residual finiteness of all finitely
presented groups.
\end{corollary}

This last corollary is in accordance with the work of Stallings \cite{stallings} which shows
that all kids of ``nasty'' subgroups can occur in $F_2\times F_2'$.

\section*{Acknowledgement}
The first author would like to thank Armando Martino for various conversations
concerning the profinite topology of $F_2\times F_2$. Both authors would also like to thank the 
anonymous referee for his careful reading and suggestions and for the simplification of the
proof of Corollary 1.

\end{document}